\documentclass[twoside, 12pt]{article}
\usepackage{amsfonts}
\usepackage{amsmath}
\usepackage{amssymb}
\usepackage{multirow}
\usepackage{graphics}
\usepackage{enumerate}

\newcommand{\R}{\mathbb{R}}

\newtheorem{theorem}{Theorem}[section]
\newtheorem{remark}[theorem]{Remark}
\newtheorem{example}[theorem]{Example}
\newtheorem{lemma}[theorem]{Lemma}
\newtheorem{corollary}[theorem]{Corollary}
\newtheorem{definition}[theorem]{Definition}
\newtheorem{proposition}[theorem]{Proposition}

\newtheorem{conjecture}{Conjecture}

\newcommand{\bea}{\begin{eqnarray}}
\newcommand{\eea}{\end{eqnarray}}
\newcommand{\beq}{\[\begin{aligned}}
\newcommand{\eeq}{\end{aligned}\]}

\newcommand{\comment}[1]{}
\def \bpm{\begin{pmatrix}}
\def \epm{\end{pmatrix}}

\def \bd{\begin{definition}}
\def \ed{\end{definition}}
\def \bcc{\begin{conjecture}}
\def \ecc{\end{conjecture}}
\def \bt{\begin{theorem}}
\def \et{\end{theorem}}
\def \bl{\begin{lemma}}
\def \el{\end{lemma}}
\def \bc{\begin{corollary}}
\def \ec{\end{corollary}}
\def\be#1\ee{\begin{align}#1\end{align}}

\def \ben{\begin{enumerate}}
\def \een{\end{enumerate}}

\def \ba{\begin{array}}
\def \ea{\end{array}}
\def \bp{\begin{proposition}}
\def \ep{\end{proposition}}
\def \bx{\begin{example}}
\def \ex{\end{example}}
\def \br{\begin{remark}}
\def \er{\end{remark}}
\def \bdsc{\begin{description}}
\def \edsc{\end{description}}
\def\pf{{\it \bf Proof. }}
\def \qed {\hfill \vrule height6pt width6pt depth0pt}
\def\hs{\hspace{.3cm}}

\def\ds{\displaystyle}
\def\1{1\!\!1}
\def\J{\mathbb{J}}
\def\x{\mathbf{x}}
\def\r{\rho}

\begin{document}

\title{A study on resistance matrix of graphs}
\author{ Deepak Sarma \\
 Department of Mathematical Sciences,\\
 Tezpur University,\\
 Tezpur-784028, India\\
 E-mail: deepaks@tezu.ernet.in\\
 Tel: +917896809545 }

\date{}

\pagestyle{myheadings} \markboth{Deepak Sarma}
{A study on resistance matrix of graphs}

\maketitle

\begin{abstract} In this article we consider resistance matrix of a connected graph. For unweighted graph we study some necessary and sufficient conditions for resistance regular graphs. Also we find some relationship between Laplacian matrix and resistance matrix in case of weighted graphs where all edge weights are positive definite matrices of given order.

\bigskip \noindent { \bf Keywords:} Connected graph,  Laplacian matrix,  Resistance distance, Moore Penrose inverse, (1)-inverse.  \\
\noindent {\bf Mathematics Subject Classification (2010):} 05C50,  05C12,   05B20.
\end{abstract}

\section{Introduction} \label{intro}

Throughout the article all our graphs are finite, undirected, connected and simple. By order of a graph we mean the number of vertices in the graph.
If $i,j \in V(G)$ in a graph $G$ then we write $i\sim j$ or $ij \in E(G)$ to mean that $i$ is adjacent to $j$ in $G.$ If $i$ is a vertex of a graph $G$ then by $N(i)$ we denote the set of all vertices in $G$ adjacent to $i$ and we call it the neighbourhood of $i.$ The degree of a vertex $i\in V(G)$ is the cardinality of $N(i)$  and is denoted by $d_i.$ A graph $G$ is said to be regular or degree regular if all the vertex degrees are equal. In that case if all the vertex degrees are equal to $d$ then we call the graph as $d- regular$ graph.  If $V_1\subseteq V(G)$ and $E_1\subseteq E(G),$ then by $G-V_1$ and $G-E_1$ we mean the graphs obtained from $G$ by deleting the vertices in $V_1$ and the edges $E_1$ respectively. In particular case when $V_1=\{u\}$ or $E_1=\{e\},$ we simply write $G-V_1$ by $G-u$ and $G-E_1$ by $G-e$ respectively. Tree is a connected graph having no cycle and a $k-forest$ is a union of $k$ disjoint trees.

By $\mathbb{M}_{m,n}$ we denote the class of all matrices of size $m\times n.$ Also by $\mathbb{M}_{n}$ we denote the class of all square matrices of order $n.$ For $M\in \mathbb{M}_{n}$ we write  $m_{ij}$ or $M_{ij}$ to denote the $ij-th$ element of $M.$ We set $[n]= \{1,2,\ldots , n \} .$ For $S\subset [n],$ we write $M(S)$ to denote the submatrix of $M$ obtained by removing rows and columns corresponding to $S.$ By eigenpair $(\lambda , \x)$ of a square matrix $M$ we mean that $\lambda$ is an eigenvalue of $M$ with corresponding eigenvector $\x .$ By an eigenvector we always mean the normalized eigenvector and when we talk about all the eigenvectors of a matrix we consider them as normalized and mutually orthogonal. By $\J $ and $\1$ we mean the matrix of all one's and vector of all one's respectively
 of suitable order. We will mention its order wherever its necessary. Adjoint of a matrix $M\in \mathbb{M}_n$ will be denoted by $adj (M).$

\bd
The adjacency matrix of a graph $G$ denoted by $A=A(G)$ is defined as
\beq
 A_{ij} =
  \begin{cases}
   1 & \text{if } (i,j) \in E \\
   0       & otherwise.
  \end{cases}
\eeq
\ed

 \bd
 The degree matrix of a graph $G$ denoted by $D=D(G)$ is a diagonal matrix such that
$D_{ii}=d_i.$\ed

\bd
The Laplacian matrix of a graph $G$ denoted by $L=L(G)$ is defined as
$L = D - A$.
\ed

\bd \label{MP}
Moore-Penrose inverse of ${\ds M\in \mathbb {M}_{m,n}}$ is defined as a matrix ${\displaystyle M^{+}\in \mathbb {M}_{n,m}}$ satisfying the following four conditions
\ben[(i)]
\item $MM^{+}M=M,$
\item $M^{+}MM^{+}=M^{+},$
\item $(MM^{+})^{*}=MM^{+}$ and
\item $(M^{+}M)^{*}=M^{+}M.$
\een
\ed

For ${\ds A\in \mathrm {M}_{m,n}}$ if the matrix $C$ satisfies (i) of definition (\ref{MP}),
 then we say $C$ to be a $( 1)-inverse$ of $A$ and we write it by $C=A^{(1)}.$

\bd The \textit{resistance distance} denoted by $r_{ij}$ between vertices $v_i$ and $v_j$  of a graph is defined as the effective resistance between the two vertices computed with Ohm's law, when the graph is viewed as an electrical network with each edge carrying unit resistance and attaching a battery at vertex $v_i$ and $v_j$. \ed

\bd
If $S$ be a set then a function $\r:S\times S \rightarrow \R$ is said to be distance function if $\r$ satisfies the following conditions for any $x,y,z \in S$
\ben[(i)]
\item $\r(x,y)\ge 0$
\item $\r(x,y)=0 \Leftrightarrow x=y$
\item  $\r(x,y)=\r(y,x)$
\item $\r(x,y)+\r(y,z)\ge \r(x,z).$
\een
\ed

 The \textit{resistance matrix} of a connected graph $G$ of order $n$ is defined to be $R(G)=(r_{ij}),$
where $r_{ij}$ is the resistance distance between the vertices $i$ and $j$ in $G.$
We call a connected graph $G$ to be resistance regular if all the row sums of $R(G)$ are equal.

 \par We now outline the contents of this article. In section \ref{rreg}, we establish some necessary and sufficient conditions for a connected graph to be resistance regular. Also we discuss some properties of resistance regular graph in that section. In section \ref{mwg} we study weighted connected graphs where all the weights are positive definite matrices of same order. We consider Laplacian matrix and resistance matrix and establish some relationships between them for matrix weighted connected graphs. Also we show that a matrix weighted tree is completely determined from its resistance matrix..

\section{Resistance regular graphs} \label{rreg}

\bl \cite{GX} \label{inverseentry}
Let $M$ be a square matrix of order $n$ and $(\lambda_k, c_k); \, k=1, \ldots , n$ are its eigenpairs, then
\beq   M_{ij}^+&= {\ds \sum_{k=1}^n } g(\lambda_k)c_{ik}c_{jk},
\text{ where } \\
g(\lambda_k)&= \begin{cases}
\frac{1}{\lambda_k} \text{ if } \lambda_k \ne 0 \\
0 ,  \text{ otherwise. }
\end{cases}\eeq \el

\bl \cite{GX} \label{eigenpair}
For any connected graph $G$ \beq r_{ij}={\ds \sum_{\lambda_k\ne 0}}\frac{1}{\lambda_k}(c_{ki}-c_{kj})^2  \eeq
where $(\lambda_k, c_k)$ is an eigenpair of $L(G).$
\el{}

\bd
The Kirchhoff index $Kf(G)$ of a connected graph $G$ is defined by \beq  Kf(G)={\ds \frac{1}{2}\sum_{i,j \in V(G)}}r_{ij}.\eeq
\ed

\bl \cite{ZWB} \label{Kf}
Let $G$ be a connected graph of order $n$ with $\lambda_k$ as eigenvalue of $L(G)$ for $k=1,\dots , n.$ Then
\beq Kf(G)=n{\ds \sum_{\lambda_k\ne 0}}\frac{1}{\lambda_k}. \eeq
\el

\bt \label{rowsum}
For a connected graph $G$ if $R_i$ denote the $i-th$ row sum of $R(G)$ and $L_{ii}^+$ denote the $i-th$ diagonal entry of $L^+(G),$ then $R_i=\frac{Kf(G)}{n}+nL_{ii}^+$ for each $i\in V(G).$
\et
\pf If $(\lambda_k, c_k)$ is an eigenpair of $L(G),$ then we have from Lemma (\ref{eigenpair})
\beq r_{ij}={\ds \sum_{\lambda_k\ne 0}}\frac{1}{\lambda_k}(c_{ki}-c_{kj})^2 . \eeq
 Therefore we get
 \beq
 R_i &={\ds \sum_{j\in V(G)}\sum_{\lambda_k\ne 0}}\frac{1}{\lambda_k} (c_{ki}-c_{kj})^2 \\
 &={\ds \sum_{\lambda_k\ne 0}\Big (\sum_{j\in V(G)}}\frac{1}{\lambda_k}(c_{ki}-c_{kj})^2\Big ) \\
 &={\ds \sum_{\lambda_k\ne 0} \bigg [ \frac{1}{\lambda_k}\bigg (nc_{ki}^2+\sum_{j\in V(G)}c_{kj}^2-2c_{ki}\sum_{j\in V(G)}c_{kj}\bigg)   \bigg] }\\
 &={\ds \sum_{\lambda_k\ne 0} \bigg [ \frac{1}{\lambda_k}\bigg (nc_{ki}^2+1-2c_{ki}\times 0 \bigg ) \bigg] }\\
 &={\ds \sum_{\lambda_k\ne 0}\frac{1}{\lambda_k} }+n{\ds \sum_{k\in V(G)}\frac{c_{ki}^2}{\lambda_k} }\\
 &=\frac{Kf(G)}{n}+nL_{ii}^+, \quad \text{ using  Lemmas } (\ref{Kf}) \text{ and } (\ref{inverseentry}).
 \eeq
\qed

As a Corollary to the above Theorem (\ref{rowsum}) we get the following result.
\bc \cite{ZWB} \label{main} A graph $G$ of order $n$ is resistance regular if and only if  $L_{11}^+=\ldots = L_{nn}^+.$
\ec

\bl ({Matrix-Tree Theorem}) \cite{Bapat} \label{MT}
Let $G$ be a graph and $W\subset V(G),$ then $det \, L(W)$ equals the number of spanning forests of $G$ with $|W|$ components in which each component contains one vertex of $W.$
\el

\bl \cite{Bapat} \label{rl}
For a connected graph $G$ with $L=L(G),$
\beq r_{ij}=L_{ii}^++L_{jj}^+-2L_{ij}^+=\frac{det\, L(i,j)}{det \, (i)} . \eeq \el

\bt
For any connected graph $G$ of order $n$ if $L_{ii}^+$ denote the $i-th$ diagonal entry of $L(G)$ and $t$ denote the number of spanning forests in $G$ having two components with $'i'$ in one the components is $t\Big( \frac{Kf(G)}{n}+nL_{ii}^+  \Big).$
\et
\pf If $R_i$ denote the $i-th$ row sum of $R(G)$ then from Lemma (\ref{rl}) we have
\beq R_i={\ds \sum_{j\in V(G)}r_{ij} }=\frac{1}{t}{\ds \sum_{j\in V(G)}L(i,j) }. \eeq
Now by Matrix-Tree Theorem $det \, L(i,j)$ represents the number of spanning forests in $G$ with two components with $'i'$ in one of them and $'j'$ on the other. Hence on using Theorem (\ref{rowsum}), we are done. \qed

\bt If $t$ be the number of spanning trees in a connected graph $G$ then the number of spanning forests in $G$ with two components is equal to $t(Kf(G)).$
\et
\pf By Matrix-Tree Theorem $det \, L(i,j)$ represents the number of spanning forests in $G$ with two components with $'i'$ in one of them and $'j'$ on the other. Therefore the required result follows from the facts that $r_{ij}=\frac{det\,  L(i,j)}{t} $ and
$ Kf(G)={\ds \frac{1}{2}\sum_{i\in V(G)}\sum_{j\in V(G)}}r_{ij}.$
\qed

\bl \cite{Bapat} \label{dis}
In a connected graph $G$ resistance distance is a distance function on the set $V(G).$ \el

\bl \cite{BG} \label{cut}
Let $G$ be a connected graph with a cut vertex $k.$ If vertices $i$ and $j$ lie in different component of $G-k$ then $r_{ij}=r_{ik}+r_{kj}.$ \el

\bt
Resistance regular graphs cannot have cut vertices.
\et
\pf
If possible suppose $G$ is a resistance regular graph with a cut vertex $v.$ Let $H_1$ is a component of $G-v$ with minimum number of vertices and $H_2=G-H_1-v.$ We consider $(v\sim)u\in H_1.$
\par Since $v$ is a cut vertex we have from Lemmas (\ref{cut}) and (\ref{dis})

\be
\label{cut1}
 &r_{ux}=r_{vx}+r_{uv} \quad \forall x\in H_2 \\
\label{cut2} &r_{xu}+r_{uv}\ge r_{xv} \quad \forall x \in H_1
\ee

Now
\beq
R_u&={\ds\sum _{x\in G}}\, r_{xu}\\
&=r_{uv}+{\ds \sum_{x\in H_1 \cup H_2} r_{xu}}\\
&\ge r_{uv}+{\ds \sum_{x\in H_1}(r_{xv}-r_{uv}) }+{\ds \sum_{x\in H_2}(r_{xv}+r_{uv})  } \quad \text{ using (\ref{cut1}) and (\ref{cut2})}\\
&= r_{uv}+\ds{ \sum_{x\in H_1 \cup H_2}r_{xu} }+ r_{uv} \, (|V(H_2)|-|V(H_1)|)\\
& \ge r_{uv}+R_v \qquad as \hs |V(H_2)|\ge |V(H_1)|\\
& > R_v, \text{ a contradiction.}
\eeq
Hence the result holds.
\qed

%\bd By spectral radius of a square matrix, we mean the largest absolute value of the eigenvalues of the matrix.\ed

%As for a resistance regular graph $G$ $\r(R_G)=\frac{2Kf(G)}{n},$ we get the following result.

%\bt If $G_1$ and $G_2$ are two resistance regular graphs of same order and $\r_1$ and $\r_2$ are respectively their resistance distance spectral radii then $\r_1 \ge \r_2 \Leftrightarrow Kf(G_1)\ge Kf(G_2)$ \et

\bt \label{Lij}
If $R=(r_{ij})$ be the resistance matrix of a connected graph $G$ with $R_i$ as the i-th row sum of $R,$ then
$L_{ij}^+(G)=\frac{R_i+R_j}{2n}-\frac{Kf(G)}{n^2}-\frac{r_{ij}}{2}$
\et
\pf From Theorem (\ref{rowsum}) we have
\beq  L_{ii}^+=\frac{R_i}{n}-\frac{Kf(G)}{n^2} \eeq

Now \beq r_{ij}&=L_{ii}^+ +L_{jj}^+ -2L_{ij}^+  \\
&=\frac{R_i+R_j}{n}-\frac{2Kf(G)}{n^2}-2L_{ij}^+ \\
\Rightarrow L_{ij}^+&=\frac{R_i+R_j}{2n}-\frac{Kf(G)}{n^2}-\frac{r_{ij}}{2}. \eeq
\qed

\bc \label{rLij}
If $G$ be a resistance regular graph, then
 $L_{ij}^+(G)=\frac{Kf(G)}{n^2}-\frac{r_{ij}}{2}.$
\ec
\pf In Theorem (\ref{Lij}), putting $R_i=R_j=\frac{2Kf(G)}{n}$ we get
\beq
L_{ij}^+&=\frac{4Kf(G)}{n\times 2n}-\frac{Kf(G)}{n^2}-\frac{r_{ij}}{2}\\
&=\frac{Kf(G)}{n^2}-\frac{r_{ij}}{2}
\eeq
\qed

In Theorem (\ref{Lij}) taking summation over the neighbourhood of any given vertex $i\in V(G)$ of a graph $G$ and using ${\ds \sum_{i\in V(G)} }R_i=2\, Kf(G)$ we get the following Corollary.
\bc{}
If $R_i$ be the $i-th$ row sum of $R(G)$ for a connected graph $G,$ then
\beq
{\ds \sum_{j\in N(i)}}L_{ij}^+=\frac{d_i}{n} \bigg [ \frac{R_i}{2}-\frac{Kf(G)}{n} \bigg]+\frac{1}{2} \bigg[ \frac{1}{n}{\ds \sum_{j\in N(i)} }R_j -{\ds \sum_{j\in N(i)} }r_{ij} \bigg] \eeq
where $N(i)$ is the set of all vertices in $G$ adjacent to $i.$ \ec{}

\bl \cite{Bapat} \label{eigL}
If $G$ be a graph of order $n$ and the eigenvalues of $L(G)$  be $\lambda_1\ge \lambda_2\ge \ldots \ge \lambda_{n-1}>\lambda_n=0,$ then for any $a\in \R $  the eigenvalues of $L+a\, \J$ are $\lambda_1, \lambda_2, \ldots, \lambda_{n-1}$ and $na.$
\el

\bl \cite{Bapat} \label{span}
If $G$ be a graph of order $n$ and the eigenvalues of $L(G)$  be $\lambda_1\ge \lambda_2\ge \ldots \ge \lambda_{n-1}>\lambda_n=0,$ then the number of spanning trees of $G$ equals $\frac{\lambda_1 \lambda_2 \ldots \lambda_{n-1}}{n}.$
\el

 \bl \cite{Bapat} \label{LXr}
 If $G$ is a connected graph with $L=L(G)$ then
 \beq L^+= \Big (L+\frac{\J}{n} \Big )^{-1}&-\frac{\J}{n}.\eeq
 \el

\bc{} If $R_i$ be the $i-th$ row sum of $R(G)$ for a connected graph $G,$ then the $ij-th$ cofactor of $(L+\frac{\J}{n})$ equals
\beq t \Big[1+\frac{R_i+R_j}{2}-\frac{Kf(G)}{n}-\frac{n\, r_{ij}}{2}\Big ],  \eeq
where $t$ is the number of spanning trees in $G.$
\ec
\pf From Lemma (\ref{LXr}) we have
\be{} \label{lj1} (L+\frac{\J}{n})^{-1}&=L^++\frac{\J}{n} \ee
Again from Lemmas (\ref{eigL}) and (\ref{span}), we get
\be  \quad \label{lj2}  det\, \Big(L+\frac{\J}{n}\Big)=t\,n
\ee{}
where $t$ is the number of spanning trees of $G.$

Now using (\ref{lj2}) in (\ref{lj1}) we get
\be \label{lj3}
\frac{adj \, (L+\frac{\J}{n})}{t\, n}= L_{ij}^++\frac{1}{n}
\ee

Applying Theorem (\ref{Lij}) in (\ref{lj3}),

\beq  det \,  \Big(L+\frac{\J }{n} \Big)(i,j)=(-1)^{i+j}t \Big[1+\frac{R_i+R_j}{2}-\frac{Kf(G)}{n}-\frac{n\, r_{ij}}{2}\Big ]. \eeq
Hence the result follows.
\qed

\bt
If $G$ is $k$ resistance regular graph of order $n,$ then $L_{ij}^+(G)=\frac{1}{2}(\frac{k}{n}-\frac{s_{ij}}{t}), $ where $t$ is the number of spanning trees in $G$ and $s_{ij}$ is the number of spanning  $2-forests$ in $G$ separating $i$ and $j.$
\et
\pf If $G$ is $k$ resistance regular graph of order $n,$ then $Kf(G)=\frac{k\, n}{2}$ and therefore from  Corollary (\ref{rLij}), we get
\be \label{sij} L_{ii}^+=\frac{Kf(G)}{n^2}=\frac{k}{2\, n} \quad \forall \, i \in [n] \ee
Now by Lemma (\ref{MT}),

\be \label{sij2} s_{ij}=det \, L(i,j).\ee

Also from Lemmas (\ref{rl}) and (\ref{MT}) for any connected graph we have $r_{ij}=\frac{det \, L(i,j)}{t}$ and from Corollary (\ref{main})  for resistance regular graph $r_{ij}=2\, (L_{ii}^+-L_{ij}^+ ).$ Comparing these we get
\beq L_{ij}^+&= L_{ii}^+-\frac{det \, L(i,j)}{2\, t}\\
&= \frac{k}{2\, n}-\frac{s_{ij}}{2\, t}, \quad \mbox{using }   (\ref{sij}) \text{ and } (\ref{sij2}).
\eeq
Hence the result follows.
\qed

\bd
A connected graph $G$ is said to be \textit{equiarboreal} if the number of spanning trees containing a given edge in $G$ is independent of the choice of edge.
\ed

\bl \cite{ZSB} \label{equi}
A connected graph $G$ is equiarboreal if and only if $r_{ij}=r_{uv}$ for any $ij, uv \in E(G).$
\el

\bl \cite{ZWB} \label{adj}
A connected graph $G$ of order $n$ with $R=R(G)$ is resistance regular if and only if
\beq {\ds \sum_{j\sim i} } r_{ij}= 2-\frac{2}{n} \quad \forall \, i\in V(G). \eeq  \el

\bt
An equiarboreal graph is resistance regular if and only if it is degree regular.
\et
\pf If $G$ be an equiarboreal graph, then let  $k=r_{ij}$ for $ij\in E(G).$
\beq
\therefore \qquad {\ds \sum_{j\sim i}}r_{ij}&=d_i\, k \eeq
Thus if $R_i$ denote the $i-th$ row sum of $R(G),$ then using Lemma (\ref{adj}) we get

\beq
  R_i&=R_l \\
  \Leftrightarrow {\ds \sum_{j\sim i}}r_{ij} &= {\ds \sum_{j\sim l}}r_{lj} \\
\Leftrightarrow \quad d_i &=d_l.
\eeq
Hence the result.
\qed{}

\bl \cite{j} \label{eigminor}
For $M\in \mathbb{M_n}$ the sum of the principal minors of $M$ of order $s\in [n]$ equals the sum of the products of the eigenvalues of $M$ taken $s$ at a time.
\el

\bl \cite{Bapat} \label{invertible}
Let $G$ be a connected graph of order $n$ and $L=L(G).$ If $M$ is any proper principal submatrix of $L,$ then $M^{-1}$ exist and is entrywise nonnegative matrix.
\el

\bl \cite{GX} \label{mix}
If $G$ be a connected graph of order $n$ and $L=L(G)$ then
\beq L^+ \J= 0=\J \, L^+, \hs LL^+={ \bf I}-\frac{\J}{n}=L^+L.\eeq

\el

\bd By bottleneck matrix at vertex $v$ of a graph $G$ with $L=L(G)$ we mean the matrix $\big(L(v)\big )^{-1}$ i.e. the inverse of the principal submatrix of $L$ obtained by removing row and column corresponding to $v.$
\ed

\bt \label{iff} For a connected graph $G$ of order $n$ if $L=L(G)$ and $(\lambda_k, c_k)$ are eigenpairs of $L$  then the following are all equivalent.
\ben[(i)]
\item $G$ is resistance regular.
\item  $L_{ii}^+(G)=\frac{Kf(G)}{n^2}$ for all $i \in V(G).$
%\item $n\, B=B\J+\J B,$ where $B$ is a bottleneck matrix at any vertex of $G.$
\item $r_{ij}+2L_{ij}^+(G)$ is constant for all $i,j \in V(G).$
\item  $R=2(L_{11}^+J-L^+).$
\item  $||B_i||_{l_1}$ is constant for all $i \in V(G),$ where $B_i$ is the bottleneck matrix of $G$ at vertex $i$ and $||M||_{l_1}$ is the sum of all entries of $M.$
\item  ${\ds \sum_{\lambda_i\ne 0} }\frac{c_{ij}^2}{\lambda_i} $ is constant for all $i \in V(G),$ where $\lambda_k$ is a eigenvalue of $L(G)$ with corresponding eigenvector $(c_{1k},\ldots,c_{nk})^T.$
\item  All principal minors of $L+\frac{J}{n}$ of order $n-1$ are equal to $t(1+\frac{Kf(G)}{n}),$  where t is the number of spanning trees in $G.$

\item  ${\ds \sum_{j\in N(i)}L_{ij}^+(G) }=\frac{d_iKf(G)}{n^2}+\frac{1}{n}-1 $ for all $i \in V(G).$
\een
\et

\pf $(i) \Leftrightarrow (ii)$ From Corollary (\ref{main}) we have $(ii)\Rightarrow (i).$ Again from the same Corollary (\ref{main}), $(i)$ implies  \beq L_{11}^+=\ldots = L_{nn}^+=l (say).\eeq
\qquad Then \beq  n\, l&={\ds \sum_{i\in V(G)}}  L_{ii}^+\\
&={\ds \sum_{\lambda_k\ne 0} }\frac{1}{\lambda_k} \text{ as eigenvalues of $L^+$ are 0  and } \frac{1}{\lambda_k} ,  \text{ for }\lambda_k\ne 0 \\
&=\frac{Kf(G)}{n} \text{ using Lemma } (\ref{Kf})  \\
\Rightarrow \quad l&= \frac{Kf(G)}{n^2}
\eeq
Thus $(i)\Rightarrow (ii).$

$(i) \Leftrightarrow (iii)$ First $(iii)$ implies

\beq
r_{ij}+2L_{ij}^+(G) &=c\, (say) \quad \forall i,j \in V(G)\\
\Rightarrow{\ds \sum_{j\in V(G)} }r_{ij}+2\, {\ds \sum_{j\in V(G)} }L_{ij}^+ &={\ds \sum_{j\in V(G)} }c \\
\Rightarrow R_i+0 &=c\, n \quad \text{ as } L^+\J=0 \text{ by Lemma (\ref{mix})  } \\
\Rightarrow (i)
\eeq
Again using Corollary (\ref{main}), $(i)$ implies
$$L_{11}^+=\ldots = L_{nn}^+=l (say).$$

\beq
\therefore \quad r_{ij}&=L_{ii}^++L_{jj}^+-2L_{ij}^+ \\
&= 2(l-L_{ij}^+)\\
\Rightarrow r_{ij}+2L_{ij}^+&=2\, l, \text{ a constant } \forall \, i,j \in V(G).
\eeq
$(iii) \Leftrightarrow (iv)$ is obvious from the proof of $(i) \Leftrightarrow (iii).$

$(iv) \Leftrightarrow (v)$ From \cite{BKL} we have
\be \label{btl1} L^+=(L+\J)^{-1}-\frac{1}{n^2} \, \J \ee

Again for a non singular matrix $M,$ we have from \cite{Bapat}
\be \label{btl2} det \, (M+\J)&=det \, M (1+\1 ^T M^{-1}\1) \ee
Therefore using (\ref{btl1}) in Corollary (\ref{main}) we get that $G$ is resistance regular if and only if all principal minors of order $n-1$ of $L+\J$ are equal.
\beq
&\Leftrightarrow (L+\J )(i,i) \text{ is constant } \forall \, i\in V(G) \\
&\Leftrightarrow t(1+||B_i||_{l_{1}}) \text{ is constant } \forall \, i\in V(G) \quad \text{ by (\ref{btl2})  and Lemma (\ref{invertible}) }\\
&\Leftrightarrow ||B_i||_{l_{1}} \text{ is constant } \forall \,i\in V(G).
\eeq

$(i) \Leftrightarrow (vi)$ Immediately follows from Corollary (\ref{main}) and Lemma (\ref{inverseentry}).

$(i) \Leftrightarrow (vii)$ If $\lambda_1\ge \lambda_2\ge \ldots \ge \lambda_{n-1}>\lambda_n=0$ are the eigenvalues of $L$ then from Lemma (\ref{eigL})  we see that eigenvalues of $L+\frac{\J}{n}$ are $\lambda_1, \lambda_2, \ldots, \lambda_{n-1}, 1.$
Again from Lemma (\ref{LXr}) we have
\be \label{minor} L^+=(L+\frac{\J}{n})^{-1}-\frac{\J}{n} \ee

Now using (\ref{minor}) in Corollary (\ref{main}) we observe that a graph is resistance regular if and only if all principal minors of $L+\frac{\J}{n}$ of order $n-1$ are equal $(k, say).$ Then using Lemma (\ref{eigminor}) we get
\beq
k\, n &={\ds \sum_{i=1}^{n-1}  }  \frac{ {\ds \prod_{i=1}^{n-1} } }{\lambda_i} +t\, n \\
&= {\ds \sum_{i=1}^{n-1} }\frac{tn}{\lambda_i} +tn  \quad \text{ by Lemma (\ref{span}) } \\
\Rightarrow k&=t\Big(1+  {\ds \sum_{i=1}^{n-1}}\frac{1}{\lambda_i} \Big )\\
&= t\Big(1+ \frac{Kf(G)}{n} \Big).
\eeq

$(ii) \Leftrightarrow (viii)$ From Lemma (\ref{mix}) we have for any graph
\beq
L\, L^+&= I-\frac{\J}{n}\\
\Rightarrow d_i\, L_{ii}^+-{\ds \sum_{j\in N(i)} }L_{ij}^+ &= 1-\frac{1}{n}
\eeq
Therefore the result holds.

\qed
%%%%%%%%%%%%%%%%%%%%%%%

From theorem (\ref{iff}) we immediately get the following results.

\bc  If $G$ is a resistance regular graph, then
\beq L_{ii}^+>L_{jk}^+ \quad \forall \, i,j,k \in V(G) \text{ with } j\ne k\eeq \ec

\bc
A $d-regular$ graph $G$ of order $n$ is resistance regular if and only if
\beq {\ds \sum_{j \in N(i)}}L_{ij}^+=\frac{d}{n^2}\, Kf(G)+\frac{1}{n}-1 \quad \forall i \in V(G) \eeq
where $N(i)$ is the set of vertices adjacent to i.
\ec

\bc
A resistance regular graph $G$ is regular if and only if ${\ds \sum_{j \in N(i)}}L_{ij}^+$ is constant for all $i \in V(G) $
where $N(i)$ is the set of vertices adjacent to i.
\ec
%%%%%%%%%%%%%%%%%%
\bt
Let $u$ and $v$ be any two vertices of a graph $G$ of order $n$ and suppose that
\beq
L(G)+\frac{\J}{n}= \bpm
d_u + 1 & k \, &\vline   & x^T \\
k & d_v + 1 \,  &\vline  & y^T \\
\hline
 x & y  \, & \vline & P
\epm
\eeq
If $P$ is invertible then $G$ is resistance regular if and only if
$d_v-d_u=y^TP^{-1}y-x^TP^{-1}x.$
\et

\pf
From Lemma (\ref{LX}) and Corollary (\ref{main}) we have
\beq G \text{ is resistance regular } &\Leftrightarrow \text{ all the principal minors of } L+\frac{\J}{n} \text{ are equal.}\\
&\Leftrightarrow (d_v+1)det P -y^T(adj P)y=(d_u+1)det P-x^T(adj P)x \\
&\Leftrightarrow (d_v-d_u)det P = y^T(adj P)-x^T(adj P)x \\
&\Leftrightarrow d_v-d_u=y^TP^{-1}y-x^TP^{-1}x.
\eeq
\qed

%%%%%%%%%%%%%%%%%%%%%%%%%%%%%%%%%%%%%%%%%%%%%%%%%

\section{Resistance matrix for matrix weighted graph} \label{mwg}

In this section we consider connected graphs with every edge associated with positive definite matrices of given order. We call such a graph to be matrix weighted graph.
If $G$ be a connected graph of order $n$ and $W_{ij}$ is the weight (positive definite matrix of order k) of the edge connecting vertices $i$ and $j$ in $G,$ then the Laplacian matrix of $G$ is defined as the $n\, k\times n\, k$ matrix (block matrix) whose $ij-th$ block is given by
\beq
L_{ij}= \begin{cases}
{\ds \sum_{j\in V(G)} } W_{ij}^{-1}  &\text{ if } \,  i=j \\
-W_{ij}^{-1} & \text{ if } i\ne j  \text{ and } i\sim j \\
0 &\text{ if } i\ne j  \text{ and } i\nsim j
\end{cases}
\eeq
Also if $L^{(1)}$ is any $(1)-inverse$ of $L(G)$ then the  $ij-th$ block of the resistance matrix $R(G)$ of the graph $G$ is defined as $R_{ij}=L_{ii}^{(1)}+L_{jj}^{(1)}- L_{ij}^{(1)} - L_{ji}^{(1)} .$

If $M$ and $N$ are two block matrices of same order and same block partition then by $M\boxtimes N $ we mean the block matrix whose $ij-th$ block is given by the matrix product $M_{ij}N_{ij},$ where $M_{ij}$ and $N_{ij}$ are $ij-th$ block of the respective matrices. Again for any two matrices $M$ and $N$ by $M\otimes N$ we mean the block matrix whose $ij-th$ block is $m_{ij}N,$ where $m_{ij}$ is the $ij-th$ entry of the matrix $M.$ If $M$ is a block matrix with all blocks of same order then by $||M||$ we mean the sum of all blocks of $M,$ i.e. $||M||= {\ds \sum M_{ij} , }$ where the summation runs over all blocks of $M.$ By $M_B^T$ we mean the block transpose of $M,$ i.e.  the $ij-th$ block of  $M_B^T$ is the $ji-th$ block of $M.$ Also we write $trace_B(M)$ to denote the sum of all diagonal blocks of $M.$

Let $X=\Big ( L+\frac{1}{n}\J_n\otimes I_k  \Big )^{-1}$ and $\widetilde{ X}$ be a block matrix with

\beq
\widetilde{X}_{ij}=\begin{cases}
X_{ii} &\text{if } i=j \\
0, &\text{otherwise. }
\end{cases}
\eeq
Then it can be easily verified that $R_{ij}=X_{ii}+X_{jj}-X_{ij}-X_{ji}$ and $R=\widetilde{X}(\J_n\otimes I_k)+(\J\otimes I_k)\widetilde{X}-X-X_B^T.$

\bt \label{wL}
For matrix weighted connected graph $G$ of order $n$ with $L=L(G)$
\beq L^+=X-\frac{1}{n}\, \J_n \otimes I_k \eeq where each weight is a positive definite matrix of order k.
\et
\pf
Since $X=\Big ( L+\frac{1}{n}\J_n\otimes I_k  \Big )^{-1}$ we get
\beq
X\Big ( L+\frac{1}{n}\J_n\otimes I_k  \Big )=I_{nk}=\Big ( L+\frac{1}{n}\J_n\otimes I_k  \Big )X \\
\Rightarrow LX=I_{nk}-\frac{1}{n} \big (\J_n\otimes I_k \big )X \text{ and } XL=I_{nk}-\frac{1}{n}X \big (\J_n\otimes I_k \big ).
\eeq

Now \beq L\Big ( X-\frac{1}{n}\J_n\otimes I_k  \Big )L&=\Big ( LX-\frac{1}{n}L(\J_n\otimes I_k)  \Big )L\\
&=\Big ( I_{nk}-\frac{1}{n}(\J_n\otimes I_k)X-{\bf 0}  \Big )L \\
&= L-\frac{1}{n}(\J_n\otimes I_k)XL \\
&=L-\frac{1}{n}(\J_n\otimes I_k) \Big (I_{nk}-\frac{1}{n}(\J_n\otimes I_k) \Big)\\
&= L-\frac{1}{n}\J_n\otimes I_k + \frac{1}{n^2}(\J_n\otimes I_k)^2 \\
&=  L-\frac{1}{n}\J_n\otimes I_k + \frac{1}{n^2}\, n (\J_n\otimes I_k) \\
&=L
\eeq
Similarly the other three conditions of Moore-Penrose inverse can also be easily established.
\qed

\bt \label{LL+}
For any connected matrix weighted graph $G$ with $L=L(G)$
\beq L\, L^+=I_{nk}-\frac{1}{n}\J_n\otimes I_k= L^+\, L=L\, X  \eeq
where each weight is a positive definite matrix of order k.
\et
\pf
From definition of the matrix $X,$ and Theorem (\ref{wL}) we have
\beq{}
LL^+&=L\Big (  X-\frac{1}{n}\J_n\otimes I_k \Big )\\
&= L\, X-\frac{1}{n} L(\J_n \otimes I_k \\
&=LX-{\bf 0} \\
&=I_{nk}-\frac{1}{n}(\J_n \otimes I_k)X \\
&= I_{nk}-\frac{1}{n}\J_n \otimes I_k \quad \text{as block column sums of }  X \text{ are all equal to } I_k.
\eeq
Similarly it can be shown that $ L^+L=I_{nk}-\frac{1}{n}\J_n\otimes I_k.$
\qed

\bd
The Kirchhoff index $Kf(G)$ of a connected matrix weighted graph $G$ is defined by $ Kf(G)={\ds \frac{1}{2}\sum_{i\in V(G)}\sum_{j\in V(G)}}R_{ij},$ where $R_{ij}$ is the $ij-th$ block of $R(G).$
\ed

\bt For matrix weighted connected graph $G$ of order $n$ with $L=L(G)$
\beq Kf(G)=n {\ds \sum_{i\in V(G)}}L_{ii}^+  \eeq
\et
\pf
By definition we have
\beq
2 Kf(G)&={\ds \sum_{i,j \in V(G)}}R_{ij}\\
&={\ds \sum_{i,j \in V(G)}}(L_{ii}^++L_{jj}^+-L_{ij}^+-L_{ji}^+ )\\
&=2(n-1){\ds \sum_{i\in V(G)}}L_{ii}^+-2\, {\ds \sum_{i\ne j} }L_{ij}^+\\
&=2(n-1){\ds \sum_{i\in V(G)}}L_{ii}^++2{\ds \sum_{i\in V(G)}}L_{ii}^+\\
&=2n {\ds \sum_{i\in V(G)}}L_{ii}^+
\eeq
Hence the result follows. \qed

\bt \label{LR}
For matrix weighted connected graph $G$ of order $n$ with $L=L(G)$ and $R=R(G)$
\beq ||L\boxtimes R||=-2(n-1)I_k \eeq where each weight is a positive definite matrix of order k.
\et
\pf
We have
\beq
\Big ( L+\frac{1}{n}\J_n\otimes I_k  \Big )X&=I_{nk}\\
\Rightarrow {\ds \sum_{j=1}^n }L_{ij}X_{ji}+\frac{1}{n}{\ds \sum{j=1}^n }X_{ji}&=I_k \quad \forall i\in V(G) \\
\Rightarrow {\ds \sum_{j=1}^n }L_{ij}X_{ji}+\frac{1}{n}I_k&=I_k \quad \forall i\in V(G)\\
\Rightarrow {\ds \sum_{j=1}^n }L_{ij}X_{ji} &=\Big ( 1-\frac{1}{n} \Big )I_k \quad \forall i\in V(G) \eeq

\be \label{LX} \therefore \quad ||L\boxtimes X ||=n\Big ( 1-\frac{1}{n} \Big ) I_k \ee

Now
\beq  R&=\widetilde{X}(\J_n\otimes I_k)+(\J\otimes I_k)\widetilde{X}-X-X_B^T\\
\Rightarrow LR &=\widetilde{X}(\J_n\otimes I_k)-(I_n\otimes I_k-\frac{1}{n}\J_n \otimes I-k )-LX_B^T \quad \text{ using Theorem } (\ref{LL+}) \text{ and } L(\J\otimes I_k)= {\bf 0}\\
\eeq
Again \beq ||L\boxtimes R ||&={\ds \sum_{i,j\in V(G)}}L_{ij}R_{ij} \\
&= trace_B(LR)\\
&= ||L\widetilde{X} ||+{\bf 0}-trace_B (I_n \otimes I_k -\frac{1}{n}\J_n \otimes I_k )-||L\boxtimes X ||\\
&= {\ds \sum_{j=1}^n }\Big ( {\ds \sum_{i=1}^n}L_{ij} \Big )X_{jj}-n\Big ( 1-\frac{1}{n} \Big ) I_k-n\Big ( 1-\frac{1}{n} \Big ) I_k \quad \text{ using } (\ref{LX})\\
&= {\bf 0}-2(n-1)I_k \eeq
Hence the result.
\qed

For matrix weighted connected graph $G$ of order $n$ with $L=L(G)$ and $R=R(G)$ let us define\\
$\tau_i =2 I_k +{\ds \sum_{j\in V(G)}} L_{ij}R_{ij} $ and $\tau =(\tau_1, \tau_2, \ldots , \tau_n )^T.$

\bt
If $G$ is a matrix weighted connected graph of order $n$ with $L=L(G)$ and $R=R(G),$ then \beq
\big (\1^T \otimes I_k   \big )\tau =2\, I_k, \text{ where each weight matrix is of order k and $\tau $ is as defined above.}
\eeq
\et
\pf We have
\beq
\tau_i &=2 I_k +{\ds \sum_{j\in V(G)}} L_{ij}R_{ij} \\
\Rightarrow {\ds \sum_{i=1}^n }\tau_i &= 2\, n I_k +||L\boxtimes R || \\
\Rightarrow  (\1^T \otimes I_k  )\tau &= 2\, nI_k -2(n-1)I_k \hs \text{ using Theorem (\ref{LR}) } \\
&= 2I_k. \eeq \qed

\bl \cite{AKB} \label{rank}
If $T$ be a matrix weighted tree on $n$ vertices such that each weight matrix is of order $k,$ then rank of $L(T)$ is $(n-1)k.$
\el

Similar to (\cite{ZWB} Theorem (9) ) we get the following for matrix weighted graph.
\bt \label{1inv}
If $G$ be a matrix weighted tree of order $n,$ then
$\bpm
L(u)^{-1} &\vline & 0 \\
\hline
0 & \vline & 0
\epm $
is a $(1 )-inverse$ of $L$, where $u$ is the vertex corresponding to the last row(column) of $L.$
\et

\pf
Let $ L= \bpm
L(u)^{-1} &\vline & \x \\
\hline
\x^T & \vline & L_{uu}
\epm  ,$ where $L_{uu}=L_{nn}.$
By Schur complement formula we have
\beq rank(L)=rank(L(u))+rank \big (L_{uu}-\x^T(L(u))^{-1} \x  \big ). \eeq

But from Lemma (\ref{rank}) we see that $rank (L)=(n-1)k=rank(L(u))$ and therefore
\be \label{sym}
L_{uu}=\x^TL(u)^{-1}\x. \ee
Now
\beq
L \bpm
L(u)^{-1} &\vline & 0 \\
\hline
0 & \vline & 0
\epm
&= \bpm{}
I_{(n-1)k} & 0 \\

\x^T L(u)^{-1} &0
\epm \\
\\
\Rightarrow \quad L \bpm
L(u)^{-1} &\vline & 0 \\
\hline
0 & \vline & 0
\epm L
&= \bpm
L(u) & \x \\
\x^T & \x^T L(u)^{-1}\x
\epm \\
&=L. \qquad \text{ using } (\ref{sym})
\eeq
Hence the result.
\qed

\bt A matrix weighted tree is completely determined from its resistance matrix. \et
\pf
 If $G$ is a matrix weighted tree then from Theorem (\ref{1inv}) we have  $\bpm
L(u)^{-1} &\vline & 0 \\
\hline
0 & \vline & 0
\epm $
as a $(1)-inverse$ of $L(G).$

Now if $R$ is known then from the relation $R_{ij}=L_{ii}^{(1)}+L_{jj}^{(1)}-L_{ij}^{(1)} -L_{ji}^{(1)} ,$ $L(u)^{-1}$ is  known. Thus $L(u)$ is determined by $R.$ \\
Again as $L\, \J_{nk}=0=\J_{nk}L,$ the structure of $G$ and edge weights of $G$ are completely known from $R$ upto isomorphism.

\qed{}

{\bf Acknowledgement:} The financial assistance for the author was provided by CSIR, India, through JRF.

\bibliographystyle{spmpsci}      % mathematics and physical sciences

\end{document}